\documentclass[a4paper,12pt]{article}

\usepackage{amsmath}
\usepackage{amssymb}
\usepackage{amsthm}
\usepackage{graphicx}
\usepackage{amscd}
\usepackage{epic, eepic}
\usepackage{url}
\usepackage{color}
\usepackage[utf8]{inputenc} 
\usepackage{enumerate}   
\usepackage{todonotes}
\usepackage{graphicx}
\usepackage{epstopdf}
\usepackage{enumitem}
\usepackage{tikz}
\usepackage[hidelinks]{hyperref}
\usepackage[top=25mm, bottom=25mm, left=28mm, right=28mm]{geometry}
\usepackage{comment}

\makeatletter
\renewenvironment{proof}[1][\proofname] {\par\pushQED{\qed}\normalfont\topsep6\p@\@plus6\p@\relax\trivlist\item[\hskip\labelsep\bfseries#1\@addpunct{.}]\ignorespaces}{\popQED\endtrivlist\@endpefalse}
\makeatother

\newtheorem{theorem}{\bf Theorem}[section]
\newtheorem{lemma}[theorem]{\bf Lemma}
\newtheorem{corollary}[theorem]{\bf Corollary}

\newtheorem{conjecture}[theorem]{\bf Conjecture}

\theoremstyle{definition}

\def\ex{\mathrm{ex}}
\def\d{\delta}
\def \cG{\mathcal{G}}

\def\blownuptheta{\theta_{3,t}\lbrack2\rbrack}
\def\e{\varepsilon}

\title{On the Tur\'an number of the blow-up of the hexagon}

\date{}

\author{Oliver Janzer\thanks{Department of Mathematics, ETH Z\"urich, Switzerland. The author is supported by an ETH Z\"urich Postdoctoral Fellowship 20-1 FEL-35.
E-mail: {\tt oliver.janzer@math.ethz.ch}}
\and 
Abhishek Methuku\thanks{School of Mathematics, University of Birmingham,
United Kingdom. The author is supported by the EPSRC, grant no. EP/S00100X/1 (A. Methuku). Email: {\tt abhishekmethuku@gmail.com}}
\and  
Zolt\'an L\'or\'ant Nagy\thanks{MTA--ELTE Geometric and Algebraic Combinatorics Research Group,
  E\"otv\"os Lor\'and University, Budapest, Hungary. The author is supported by the Hungarian Research Grant (NKFI) No. K 120154 and by the János Bolyai Scholarship of the Hungarian Academy of Sciences. 	E-mail: {\tt nagyzoli@cs.elte.hu}}}
  
\begin{document}

\maketitle

\begin{abstract}
    The $r$-blowup of a graph $F$, denoted by $F[r]$, is the graph obtained by replacing the vertices and edges of $F$ with independent sets of size $r$ and copies of $K_{r,r}$, respectively.
    For bipartite graphs $F$, very little is known about the order of magnitude of the Tur\'an number of $F[r]$. \\
    \indent In this paper we prove that $\ex(n,C_6[2])=O(n^{5/3})$ and, more generally, for any positive integer $t$, $\ex(n,\theta_{3,t}[2])=O(n^{5/3})$. This is tight when $t$ is sufficiently large.
\end{abstract}

\section{Introduction}

The Tur\'an number (or extremal number) of a graph $H$, denoted $\ex(n,H)$, is the maximum possible number of edges in an $H$-free graph on $n$ vertices. The famous Erd\H os--Stone--Simonovits theorem \cite{ES46,ESi66} states that $\ex(n,H)=(1-\frac{1}{\chi(H)-1}+o(1)){\frac{n^2}{2}}$, which determines the asymptotics of $\ex(n,H)$ when $H$ has chromatic number at least $3$. However, for bipartite graphs, the theorem only gives $\ex(n,H)=o(n^2)$. Finding better bounds in this case is a central and difficult area in extremal graph theory. For example, a classical result of Bondy and Simonovits \cite{BS74} states that $\ex(n,C_{2k})=O(n^{1+1/k})$, but this is only known to be tight for $k\in \{2,3,5\}$. Recently, good estimates have been obtained for a wide range of graphs $H$, especially for subdivisions.

The topic of this paper is the Tur\'an number of blow-up graphs. For a graph $F$ and a positive integer $r$, the $r$-blowup of $F$ is the graph obtained by replacing the vertices and edges of $F$ with independent sets of size $r$ and copies of $K_{r,r}$, respectively. We denote this graph by $F[r]$. The study of the Tur\'an number of blow-ups was initiated by Grzesik, Janzer and Nagy \cite{GJN19}. They proved, among other things, that when $T$ is a tree, then $\ex(n,T[r])=O(n^{2-1/r})$. They have also made the following general conjecture.

\begin{conjecture}[Grzesik--Janzer--Nagy \cite{GJN19}]
    Let $F$ be a graph such that $\ex(n,F)=O(n^{2-\alpha})$ for some constant $0\leq \alpha\leq 1$. Then for any positive integer $r$,
    $$\ex(n,F[r])=O(n^{2-\frac{\alpha}{r}}).$$
\end{conjecture}

Their result mentioned above proves this conjecture when $F$ is a tree. The conjecture also holds when $F=K_{s,t}$ and $\alpha=\frac{1}{s}$. This covers the case $C_4=K_{2,2}$, but already when $F=C_{2k}$ for some $k\geq 3$, the conjecture is open. In this case, it can be stated as follows.

\begin{conjecture} \label{con:blownupcycles}
    For any $r,k\geq 2$,
    $$\ex(n,C_{2k}[r])=O(n^{2-\frac{1}{r}+\frac{1}{rk}}).$$
\end{conjecture}

\sloppy Since $C_{2k}[r]$ has maximum degree $2r$, a result of F\"uredi \cite{Fu91} (reproved by Alon, Krivelevich and  Sudakov \cite{AKS03}) gives $\ex(n,C_{2k}[r])=O(n^{2-\frac{1}{2r}})$. Note also that $C_{2k}[r]$ does not contain $K_{2r,2r}$ as a subgraph when $k>2$, thus a result of Sudakov and Tomon \cite{ST20} shows that in this case $\ex(n,C_{2k}[r])=o(n^{2-\frac{1}{2r}})$. Very recently, the first author proved \cite{Jan20} that $\ex(n,C_{2k}[r])=O(n^{2-\frac{1}{r}+\frac{1}{r+k-1}}(\log n)^{\frac{4k}{r(r+k-1)}})$, which improves this bound when $k>r+1$.

In this paper we establish Conjecture \ref{con:blownupcycles} in the first unknown case -- the $2$-blowup of the hexagon.

\begin{theorem}
    $$\ex(n,C_6[2])=O(n^{5/3}).$$
\end{theorem}

In fact, we can prove a more general result about theta graphs. The theta graph $\theta_{k,t}$ is the union of $t$ paths of length $k$ which share the same endpoints but are pairwise internally vertex-disjoint. Note that $C_{2k}=\theta_{k,2}$. Our more general result is as follows.

\begin{theorem} \label{thm:turantheta}
    For any positive integer $t$,
    $$\ex(n,\theta_{3,t}[2])=O(n^{5/3}).$$
\end{theorem}

This is tight for sufficiently large $t$ by a general result of Bukh and Conlon \cite{BC18}, which we will now state. Let $F$ be a graph with a set $R\subsetneq V(F)$ of distinguished vertices, called the roots of $F$. For any non-empty subset $S\subset V(F)\setminus R$, define $\rho_F(S)$ to be $\frac{e_S}{|S|}$, where $e_S$ is the number of edges in $F$ with at least one endpoint in $S$. Let $\rho(F)=\rho_F(V(F)\setminus R)$. We say that $F$ is balanced if $\rho(F)\leq \rho_F(S)$ for every non-empty $S\subsetneq V(F)\setminus R$. For a positive integer $t$, we define the rooted $t$-blowup of $F$ to be the graph obtained by taking $t$ pairwise vertex-disjoint copies of $F$ and, for every $v\in R$, identifying the various copies of $v$. We note that, somewhat confusingly, the notions `rooted $t$-blowup' and `$t$-blowup' are completely different. Observe that $\theta_{k,t}$ is the rooted $t$-blowup of the path of length $k$ whose roots are the leaves. Similarly, $\theta_{k,t}[2]$ is the rooted $t$-blowup of $P_k[2]$ whose roots are the $4$ vertices of degree $2$ (here and below $P_k$ denotes the path with $k$ edges). Let us write $t\ast F$ for the rooted $t$-blowup of $F$.
The result of Bukh and Conlon is as follows.

\begin{theorem}[Bukh--Conlon {\cite[Lemma 1.2]{BC18}}]
    Let $F$ be a balanced rooted graph. Then there exists a positive integer $t_0=t_0(F)$ such that for every $t\geq t_0$,
    $$\ex(n,t\ast F)=\Omega(n^{2-\frac{1}{\rho(F)}}).$$
\end{theorem}

They stated their result in the case where $F$ is a tree, but as Kang, Kim and Liu observed \cite{KKL18}, this assumption is not used in their proof.

Let $F=P_3[2]$ and let the roots of $F$ be its degree $2$ vertices. Then note that $F$ is balanced with $\rho(F)=\frac{12}{4}=3$. Moreover, as we have already remarked, $\theta_{3,t}[2]$ is the rooted $t$-blowup of this graph. Thus, it follows that for sufficiently large $t$ we have $\ex(n,\theta_{3,t}[2])=\Omega(n^{5/3})$. Together with Theorem \ref{thm:turantheta}, we get the following result.

\begin{corollary}
    For sufficiently large $t$, we have
    $$\ex(n,\theta_{3,t}[2])=\Theta(n^{5/3}).$$
\end{corollary}

\subsection{Outline of the proof}
\label{outline}

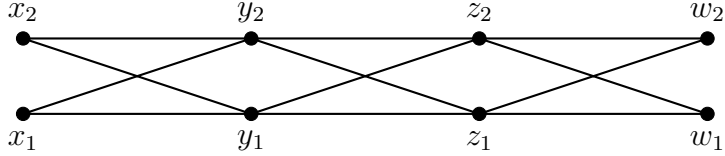
\begin{figure}
		\centering
		\begin{tikzpicture}[scale=0.5]
		\draw[fill=black](0,0)circle(5pt);
		\draw[fill=black](0,2)circle(5pt);		
		\draw[fill=black](6,0)circle(5pt);
		\draw[fill=black](6,2)circle(5pt);
		\draw[fill=black](12,0)circle(5pt);
		\draw[fill=black](12,2)circle(5pt);
		\draw[fill=black](18,0)circle(5pt);
		\draw[fill=black](18,2)circle(5pt);
		
		\draw[thick](0,0)--(6,0)--(0,2)--(6,2)--(0,0);
		\draw[thick](6,0)--(12,0)--(6,2)--(12,2)--(6,0);
		\draw[thick](12,0)--(18,0)--(12,2)--(18,2)--(12,0);
		
		\node at (0,-0.7) {$x_1$};
		\node at (0,2.7) {$x_2$};
		\node at (6,-0.7) {$y_1$};
		\node at (6,2.7) {$y_2$};
		\node at (12,-0.7)  {$z_1$};
		\node at (12,2.7)  {$z_2$};
		\node at (18,-0.7) {$w_1$};
		\node at (18,2.7) {$w_2$};
		\end{tikzpicture}
		\caption{A copy of $P_3[2]$}
		\label{fig:blownuppath}
\end{figure}

Before we get on with the proof of Theorem \ref{thm:turantheta}, let us give a brief sketch of the argument. First, using a standard reduction lemma, we will assume that our host graph $G$ is nearly regular. Then we will find many copies of $P_3[2]$ in $G$ with a fixed pair of endpoints $(x_1,x_2)$. Here and below, $P_3[2]$ has vertices $x_1,x_2,y_1,y_2,z_1,z_2,w_1,w_2$ and edges $x_iy_j,y_iz_j,z_iw_j$ (see Figure \ref{fig:blownuppath}).

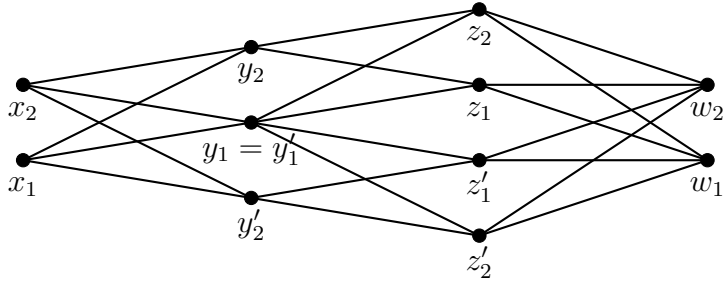
\begin{figure}
		\centering
		\begin{tikzpicture}[scale=0.5]
		\draw[fill=black](0,0)circle(5pt);
		\draw[fill=black](0,2)circle(5pt);		
		\draw[fill=black](6,-1)circle(5pt);
		\draw[fill=black](6,1)circle(5pt);
		\draw[fill=black](6,3)circle(5pt);
		\draw[fill=black](18,0)circle(5pt);
		\draw[fill=black](18,2)circle(5pt);
		\draw[fill=black](12,4)circle(5pt);
		\draw[fill=black](12,2)circle(5pt);
		\draw[fill=black](12,0)circle(5pt);
		\draw[fill=black](12,-2)circle(5pt);
		
		\draw[thick](0,2)--(6,3)--(12,4)--(18,2);
		\draw[thick](0,2)--(6,1)--(12,4)--(18,0);
		\draw[thick](0,0)--(6,3)--(12,2)--(18,2);
		\draw[thick](0,0)--(6,1)--(12,2)--(18,0);
		
		\draw[thick](0,0)--(6,-1)--(12,0)--(18,0);
		\draw[thick](0,2)--(6,-1)--(12,-2)--(18,0);
		\draw[thick](6,1)--(12,0)--(18,2);
		\draw[thick](6,1)--(12,-2)--(18,2);
		
		\node at (0,-0.7) {$x_1$};
		\node at (0,1.3) {$x_2$};
		\node at (6,-1.7) {$y_2'$};
		\node at (6,2.3) {$y_2$};
		\node at (6,0.3)  {$y_1=y_1'$};
		\node at (12,3.3)  {$z_2$};
		\node at (12,-2.7)  {$z_2'$};
		\node at (12,-0.7)  {$z_1'$};
		\node at (12,1.3) {$z_1$};
		\node at (18,-0.7) {$w_1$};
		\node at (18,1.3) {$w_2$};
		\end{tikzpicture}
		\caption{A degenerate copy of $C_6[2]$}
		\label{fig:degenerate blownupc6}
\end{figure}

It is not hard to see that if $G$ has minimum degree at least $Cn^{2/3}$ for some big constant $C$, then for some pair $(x_1,x_2)\in V(G)^2$, much more than $n^2$ such copies can be found. This means that there will be many among these copies that share the same $(w_1,w_2)$.
If we take $t$ such $P_3[2]$'s, their union is a homomorphic copy of $\theta_{3,t}[2]$. However, it may be a degenerate one, i.e. some of the internal vertices may coincide in the $t$ copies of $P_3[2]$. A simple counting argument shows that unless $G$ contains $\blownuptheta$, there must be plenty of degenerate copies of $C_6[2]$ in it. (See Figure~\ref{fig:degenerate blownupc6} for one of the possibilities. We have two copies of $P_3[2]$ there, one with vertices $x_1,x_2,y_1,y_2,z_1,z_2,w_1,w_2$ and one with vertices $x_1,x_2,y_1',y_2',z_1',z_2',w_1,w_2$, but they do not give a genuine $C_6[2]$ because $y_1=y_1'$.)

Ideally, we would like to efficiently bound the number of these degenerate copies from above and thereby reach a contradiction. However, this seems impossible without further assumptions on the ``building blocks" -- the $P_3[2]$'s. Accordingly, we will only use those $P_3[2]$'s in the above argument which satisfy some extra properties. For example, we only count those $P_3[2]$'s for which $d(x_1,x_2,z_1,z_2)<6t$ and $d(y_1,y_2,w_1,w_2)<6t$ hold for the common neigbourhoods of the quadruples in consideration. Lemma \ref{lemma:4heavy} (from Section~\ref{sec:2}) will show that we do not lose too many $P_3[2]$'s by doing so. We will also make sure that in all our $P_3[2]$'s, the codegree $d(z_1,z_2)$ is roughly the same (up to a factor of 2). Finally, we will insist that $d(x_1,x_2,z_1)$ and $d(x_1,x_2,z_2)$ are not too large (relative to $d(x_1,x_2)$). In Lemma \ref{lemma:goodblownuppaths}, we show that we have many $P_3[2]$'s possessing these properties.

\begin{figure}
		\centering
		\begin{tikzpicture}[scale=0.5]
		\draw[fill=black](0,0)circle(5pt);
		\draw[fill=black](0,2)circle(5pt);		
		\draw[fill=black](6,-1)circle(5pt);
		\draw[fill=black](6,1)circle(5pt);
		\draw[fill=black](6,3)circle(5pt);
		\draw[fill=black](18,0)circle(5pt);
		\draw[fill=black](18,2)circle(5pt);
		\draw[fill=black](12,4)circle(5pt);
		\draw[fill=black](12,2)circle(5pt);
		\draw[fill=black](12,0)circle(5pt);
		\draw[fill=black](12,-2)circle(5pt);
		
		\draw[fill=black](18,4)circle(5pt);
		\draw[fill=black](18,6)circle(5pt);
		\draw[fill=black](12,8)circle(5pt);
		\draw[fill=black](12,6)circle(5pt);
		
		\draw[thick](0,2)--(6,3);
		\draw[thick](0,2)--(6,1);
		\draw[thick](0,0)--(6,3);
		\draw[thick](0,0)--(6,1);
		
		\draw[thick](0,0)--(6,-1)--(12,0);
		\draw[thick](0,2)--(6,-1)--(12,-2);
		\draw[thick](6,1)--(12,0);
		\draw[thick](6,1)--(12,-2);
		
		\draw[green,thick](6,1)--(12,4)--(6,3)--(12,2)--(6,1);
		\draw[green,thick](12,2)--(18,0)--(12,4)--(18,2)--(12,2);
		\draw[green,thick](12,-2)--(18,0)--(12,0)--(18,2)--(12,-2);
		
		\draw[red,thick](6,1)--(12,8)--(6,3)--(12,6)--(6,1);
		\draw[red,thick](12,6)--(18,4)--(12,8)--(18,6)--(12,6);
		\draw[red,thick](12,-2)--(18,4)--(12,0)--(18,6)--(12,-2);
		
		\node at (0,-0.7) {$x_1$};
		\node at (0,1.3) {$x_2$};
		\node at (6,-1.7) {$y_2'$};
		\node at (6,2.3) {$y_2$};
		\node at (6,0.3)  {$y_1=y_1'$};
		\node at (12,3.3)  {$z_2^1$};
		\node at (12,-2.7)  {$z_2'$};
		\node at (12,-0.7)  {$z_1'$};
		\node at (12,1.3) {$z_1^1$};
		\node at (18,-0.7) {$w_1^1$};
		\node at (18,1.3) {$w_2^1$};
		
		\node at (12,5.3) {$z_1^2$};
		\node at (12,7.3)  {$z_2^2$};
		\node at (18,3.3) {$w_1^2$};
		\node at (18,5.3) {$w_2^2$};
		
		\end{tikzpicture}
		\caption{The red and green edges together give a copy of $C_6[2]$}
		\label{fig:change base}
\end{figure}
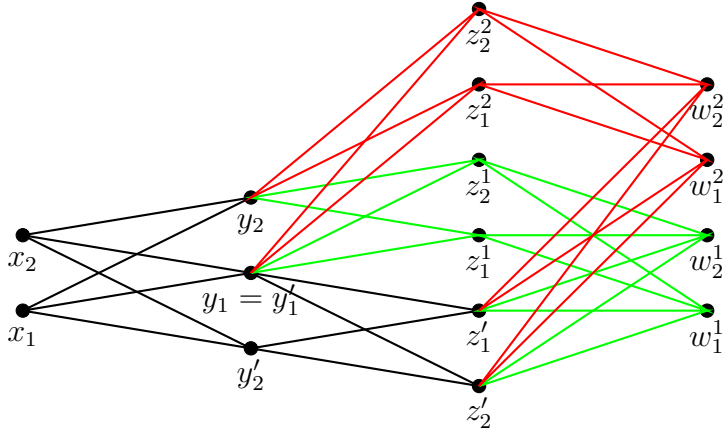

Thus we can now turn to bounding the number of degenerate copies of $C_6[2]$ which arise from these ``controlled" $P_3[2]$'s. Let us discuss briefly how we will upper bound the number of configurations depicted in Figure \ref{fig:degenerate blownupc6}. We first argue that if there are many such configurations, then there exists some choice of $y_1,y_2,y_2',z_1',z_2'$ which can be extended to many such bad configurations. However, if these vertices extend to at least $t$ such configurations and in those extensions the vertices $z_1,z_2,w_1,w_2$ are all distinct, then we get a $\blownuptheta$ in $G$ (the case $t=2$ is illustrated in Figure \ref{fig:change base}). This means that the vertices $z_1,z_2,w_1,w_2$ cannot all be distinct in the various extensions. We will argue that in fact $z_1$ (or $z_2)$ will be the same vertex in many of these extensions. This will mean that the codegree $d(z_1,z_1',z_2')$ is very large (larger than $n^{1/6}$) in these configurations.

\begin{figure}
		\centering
		\begin{tikzpicture}[scale=0.7]
		
		\draw[fill=black](3,6)circle(5pt);
		\draw[fill=black](2,6)circle(5pt);
		\draw[fill=black](-2,6)circle(5pt);
		\draw[fill=black](-5,6)circle(5pt);
		\draw[fill=black](-1,6)circle(5pt);
		\draw[fill=black](-4,6)circle(5pt);
		
		\draw[fill=black](-6,12)circle(5pt);
		\draw[fill=black](-5,12)circle(5pt);
		\draw[fill=black](-2,12)circle(5pt);
		\draw[fill=black](-1,12)circle(5pt);
		\draw[fill=black](1,12)circle(5pt);
		\draw[fill=black](2,12)circle(5pt);
					
		\draw[thick](-6,12)--(-5,6)--(-2,12);
		\draw[thick](-6,12)--(-4,6)--(-2,12);
		\draw[thick](-5,12)--(-5,6)--(-1,12);
		\draw[thick](-5,12)--(-4,6)--(-1,12);
		
		\draw[thick](-6,12)--(-2,6)--(1,12);
		\draw[thick](-6,12)--(-1,6)--(1,12);
		\draw[thick](-5,12)--(-2,6)--(2,12);
		\draw[thick](-5,12)--(-1,6)--(2,12);
		
		\draw[dashed,thick](2,6)--(1,12)--(3,6)--(2,12)--(2,6);
		\draw[dashed,thick](2,6)--(-2,12)--(3,6)--(-1,12)--(2,6);
		
		\node at (2,5.5) {$z_1'$};
		\node at (3,5.5) {$z_2'$};
		
		\draw[rotate around={90:(-3,6)}] (-3,6) ellipse (1 and 4);
		\draw[rotate around={90:(-2,12)}] (-2,12) ellipse (1 and 6);
		\node at (-10,12)  {$N(z_1',z_2')$};
		\node at (-8,6)  {$Y$};
		\end{tikzpicture}		
		\caption{Finding a $C_6[2]$ using the $P_4[2]$ formed by the solid edges}
		\label{fig:3heavy illustration}
	\end{figure}
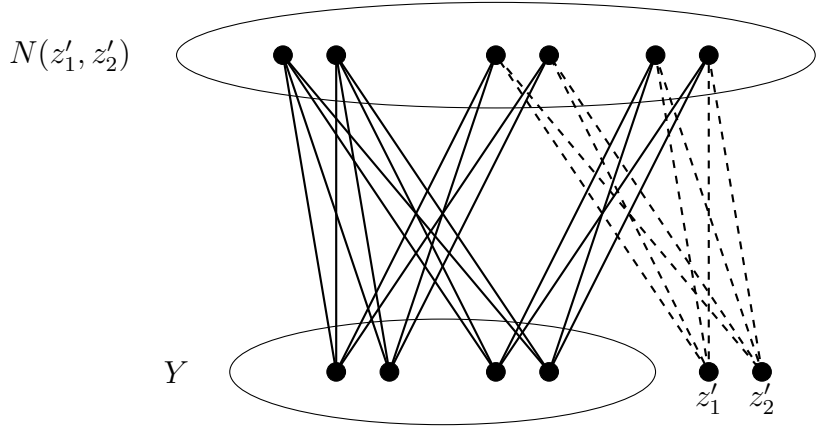

Here comes one of the key ideas of the proof (the precise statement is in Lemma~\ref{lemma:3heavy}). We give an upper bound for the number of configurations depicted in Figure \ref{fig:degenerate blownupc6} for which $d(z_1,z_1',z_2')$ is much bigger than $n^{1/6}$ roughly as follows. First note that since our graph is almost regular with edge density around $n^{-1/3}$, it is reasonable to assume that the codegree of a typical pair of vertices is around $n^{1/3}$. We claim that for a fixed pair $(z_1',z_2')$ of vertices with $d(z_1',z_2')\approx n^{1/3}$, the number of $z_1\in V(G)$ with $d(z_1,z_1',z_2')\gg n^{1/6}$ is at most $n^{1/3}$. Indeed, otherwise we can find a bipartite subgraph $B$ in $G$ with parts $X=N(z_1',z_2')$ and $Y$ where $|X|\approx |Y|\approx n^{1/3}$ and the degree of each vertex of $Y$ in $B$ is much more than $n^{1/6}$, so $B$ has much more than $n^{1/2} = |V(B)|^{3/2}$ edges. By a result of Grzesik, Janzer and Nagy \cite{GJN19}, this implies that we can find the $2$-blowup of any tree in $B$. In particular, we can find the $2$-blowup of a path of length $4$ whose degree $2$ vertices are embedded to $N(z_1',z_2')$. It is not hard to see that this means that $G$ contains $C_6[2]$ as a subgraph (see Figure~\ref{fig:3heavy illustration}). In the general case when the forbidden subgraph is $\blownuptheta$, we just need to embed a larger $2$-blownup tree into the bipartite graph.

The next section is split into four subsections. In Subsection \ref{subsec:quadruples}, we prove that in a $\blownuptheta$-free graph there cannot be locally too many quadruples with large codegree. Similarly, in Subsection \ref{subsec:triples}, we prove that in a $\blownuptheta$-free graph there cannot be locally too many triples with large codegree. In Subsection \ref{subsec:prescribeproperties}, we show that we have many $P_3[2]$'s with some prescribed properties. In Subsection \ref{subsec:countbadpairsofP3}, we first show that this implies that there are many copies of the configuration in Figure \ref{fig:degenerate blownupc6}. We then arrive at a contradiction by proving an upper bound for the number of these configurations.

\section{The proof of Theorem \ref{thm:turantheta}}
\label{sec:2}
	
\textbf{Notation.} Given a graph $G$, we write $N(v_1,\dots,v_k)$ for the {\em common neighbourhood} of vertices $v_1,\dots,v_k \in V(G)$, and we let $d(v_1,\dots,v_k)=|N(v_1,\dots,v_k)|$. Moreover, for a set $S\subset V(G)$, we write $d_S(v)=|N(v)\cap S|$. A graph $G$ is called {\em $K$-almost-regular} if $\max_{v \in V(G)} d(v) \le K\cdot \min_{v \in V (G)} d(v)$, where $d(v)$ denotes the degree of $v$.

To compare the order of magnitude of two positive functions $f(n)$ and $g(n)$, we use the following notations. We write $f(n)=O(g(n))$ if there exists a positive constant $C$ such that $\frac{f(n)}{g(n)}\leq C$ for all $n\in \mathbb{N}$.
Similarly, $f(n)=\Omega(g(n))$ means that there exists a positive constant $C$ such that $\frac{f(n)}{g(n)}\geq C$ for all $n\in \mathbb{N}$.  We use $f(n)=\Theta(g(n))$ to express that both $f(n)=O(g(n))$ and  $f(n)=\Omega(g(n))$ hold. Finally, $f(n)=o(g(n))$ and $f(n)=\omega(g(n))$ denote that $\lim_{n\rightarrow \infty}\frac{f(n)}{g(n)}=0$ and $\lim_{n\rightarrow \infty}\frac{f(n)}{g(n)}=\infty$ hold, respectively.

In what follows, we assume that $t$ is a fixed positive integer. Accordingly, the constant $C$ in the asymptotic notations $O$ and $\Omega$ can depend on $t$.

\bigskip

We start with a lemma that allows us to restrict our attention to almost regular graphs. A version of this lemma was first proved by Erd\H os and Simonovits \cite{ESi70}. We use a slight variant due to Jiang and Seiver.

\begin{lemma}[Jiang--Seiver \cite{JS12}] \label{lemmaJS}
	Let $\e,c$ be positive reals, where $\e<1$ and $c\geq 1$. Let $n$ be a positive integer that is sufficiently large as a function of $\e$. Let $G$ be a graph on $n$ vertices with $e(G)\geq cn^{1+\e}$. Then $G$ contains a $K$-almost-regular subgraph $G'$ on $m\geq n^{\frac{\e-\e^2}{2+2\e}}$ vertices such that $e(G')\geq \frac{2c}{5}m^{1+\e}$ and $K=20\cdot 2^{\frac{1}{\e^2}+1}$.
\end{lemma}

It is well known that any graph with $e$ edges contains a bipartite subgraph with at least $e/2$ edges. This observation, combined with Lemma \ref{lemmaJS}, reduces Theorem \ref{thm:turantheta} to the following statement.

\begin{theorem} \label{thm:turanthetareduced}
	Let $K$ be a constant and let $G$ be a $K$-almost-regular bipartite graph on $n$ vertices with minimum degree $\d=\omega(n^{2/3})$. Then, for $n$ sufficiently large, $G$ contains a copy of $\blownuptheta$.
\end{theorem}

\subsection{Bounding the number of quadruples with large common neighbourhood} \label{subsec:quadruples}

As mentioned in the proof outline (Section \ref{outline}), we want to find many copies of $P_3[2]$ with some extra properties. Our first prescribed property is that the quadruples $(x_1,x_2,z_1,z_2)$ and $(y_1,y_2,w_1,w_2)$ should have few common neighbours (see Figure \ref{fig:blownuppath} for the position of these vertices). The next lemma will be used to achieve this.

\begin{lemma} \label{lemma:4heavy}
	Let $G$ be a $\blownuptheta$-free graph. Let $x$, $x'$, $y$ and $y'$ be distinct vertices in $G$ and let $R\subset N(y,y')\setminus \{x,x'\}$. Then the number of pairs of distinct vertices $(z,z')\in R^2$ with $d(x,x',z,z')\geq 6t$ is at most $4t|R|$.
\end{lemma}

\begin{proof}
	Take a maximal set of pairs $(z_1,z_1'),\dots,(z_s,z_s')\in R^2$ such that $z_1,z_1',\dots,z_s,z_s'$ are all distinct and $d(x,x',z_i,z_i')\geq 6t$ for every $i$. If $s\geq t$, then we may choose $w_1,w_1',\dots,w_t,w_t'\in V(G)$ such that $x,x',y,y'$, $z_i,z_i'$ ($1\leq i\leq t$) and $w_j,w_j'$ ($1\leq j\leq t$) are all distinct, and $w_i,w_i'\in N(x,x',z_i,z_i')$ for all $1\leq i\leq t$. Then the vertices $x,x',y,y'$, $z_i,z_i'$ ($1\leq i\leq t$) and $w_j,w_j'$ ($1\leq j\leq t$) form a copy of $\theta_{3,t}[2]$, which is a contradiction.
	
	Thus, $s< t$. By maximality, for any $(z,z')\in R^2$ with $d(x,x',z,z')\geq 6t$ we have $\{z,z'\} \cap \{z_1,z_1',\dots,z_s,z_s'\}\neq \emptyset$. This leaves at most $2\cdot 2s\cdot |R|< 4t|R|$ possibilities for such $(z,z')$. 
\end{proof}

\subsection{Bounding the number of triples with large common neighbourhood} \label{subsec:triples}

We now turn to bounding the number of triples which have very large common neighbourhood. As mentioned in the introduction, this is one of the key ideas of the proof. ``Large" common neighbourhood here means that the size is something like the square root of the codegree of a typical pair of vertices, which will be around $n^{1/6}$ for our graph with density roughly $n^{-1/3}$.

\begin{lemma} \label{lemma:3heavy}
	There exists a constant $C=C(t)$ with the following property. Let $G$ be a $\blownuptheta$-free bipartite graph. Let $z_1,z_2$ be distinct vertices in $G$ and let $N(z_1,z_2)$ have size $\ell$. Let $R= \{v\in V(G)\setminus \{z_1,z_2\}: d(v,z_1,z_2)\geq C\ell^{1/2}\}$. Then the number of triples $(z',w_1,w_2)$ of distinct vertices with $z'\in R$, $w_1,w_2\in N(z',z_1,z_2)$ is at most $C\ell^2$.
\end{lemma}

To prove this lemma, we need some preliminaries.

\begin{lemma} \label{lemma:blownuptree}
	Let $T$ be a tree and let $v\in V(T)$. Then there exists a constant $C=C(T)$ with the following property. Let $G$ be a bipartite graph with parts $X$ and $Y$ of size at most $n$. Assume that $G$ has at least $Cn^2$ copies of $K_{2,2}$. Then $G$ contains a copy of $T\lbrack 2\rbrack$ with the two images of $v$ embedded in $X$.
\end{lemma}

The proof of this lemma is similar to the proof of Theorem 1.6 in \cite{GJN19}, so it is only sketched here.

\begin{proof}[Sketch of Proof.]
    Define an auxiliary bipartite graph $\cG$ whose parts are $X^{(2)}$ and $Y^{(2)}$ (the set of $2$-element subsets of $X$ and $Y$, respectively) and in which there is an edge between $\{x_1,x_2\}\in X^{(2)}$ and $\{y_1,y_2\}\in Y^{(2)}$ if and only if $x_1$, $x_2$, $y_1$ and $y_2$ form a $K_{2,2}$ in $G$. The number of vertices in $\cG$ is $\binom{|X|}{2}+\binom{|Y|}{2}\leq 2\binom{n}{2}\leq n^2$. On the other hand, the number of edges in $\cG$ is equal to the number of $K_{2,2}$'s in $G$, which is at least $Cn^2$ by assumption. Now we construct a (not uniformly) random mapping $f:V(T)\rightarrow V(\cG)$ in a way that edges of $T$ are mapped to edges of $\cG$. This gives an embedding of $T[2]$ into $G$ provided that the images of different vertices of $T$ are mapped to vertices of $\cG$ which are disjoint as subsets of $V(G)$.
    
    We map the vertices of $T$ to $V(\cG)$ one by one, starting at $v$ and in a way that when it comes to embedding a vertex $u\in V(T)$, $u\neq v$, we have embedded precisely one neighbour of $u$ so far. For every $x\in V(\cG)$, let $f(v)=x$ with probability $\frac{d_{\cG}(x)}{2e(\cG)}$. It is easy to see that with probability $1/2$, $f(v)$ is in $X^{(2)}$. When it comes to mapping some $u\neq v$ to $V(\cG)$, let $w$ be the unique neighbour of $u$ in $T$ which has already been mapped to $V(\cG)$ and let $y=f(w)$. Then we let $f(u)$ be a uniformly random neighbour of $y$ in $\cG$.
    
    It is clear that this gives a graph homomorphism $f$ from $T$ to $\cG$. We only need to check that with probability more than $1/2$, the images of different vertices of $T$ are pairwise disjoint sets. The crucial thing to notice is that for every $u\in V(T)$, the distribution of $f(u)$ is the same as the distribution of $f(v)$, namely $\mathbb{P}(f(u)=x)=\frac{d_{\cG}(x)}{2e(\cG)}$, meaning that vertices in $\cG$ whose degree is much lower than the average degree of $\cG$ are unlikely to be in the image of $f$. Using this and the fact that the average degree of $\cG$ is at least $C$, we can conclude that the probability that there exists a pair of distinct vertices in $T$ whose images are not disjoint sets in $V(G)$ tends to $0$ as $C$ tends to infinity. Thus, for $C$ sufficiently large, we get a suitable embedding with positive probability.
\end{proof}

\begin{lemma} \label{lemma:count cherries}
	Let $T$ be a tree and let $v\in V(T)$. Then there exists a constant $C=C(T)$ with the following property. Let $H$ be a bipartite graph with parts $X$ and $Y$ such that $d(y)\geq C|X|^{1/2}$ for every $y\in Y$. Assume that the number of triples of distinct vertices $(y,x_1,x_2)\in Y\times X\times X$ with $yx_1,yx_2\in E(H)$ is more than $C|X|^2$. Then $H$ contains a copy of $T[2]$ with the two images of $v$ embedded in $X$.
\end{lemma}

\begin{proof}
    Let $C$ be sufficiently large. Assume first that $|Y|\leq |X|$. Since the number of triples of distinct vertices $(y,x_1,x_2)\in Y\times X\times X$ with $yx_1,yx_2\in E(H)$ is at least $C|X|^2$, on average a pair of distinct vertices in $X$ have at least $C$ common neighbours. Hence, the number of $K_{2,2}$'s in $H$ is at least $\binom{|X|}{2}\binom{C}{2}\geq \frac{1}{1000}C^2|X|^2$. Since $|Y|\leq |X|$ and $C$ is sufficiently large, Lemma \ref{lemma:blownuptree} shows that $H$ contains a copy of $T[2]$ with the two images of $v$ embedded in $X$.
    
    Assume now that $|X|<|Y|$. Since $d(y)\geq C|X|^{1/2}$ for every $y\in Y$, the number of triples of distinct vertices $(y,x_1,x_2)\in Y\times X\times X$ with $yx_1,yx_2\in E(H)$ is at least $|Y|C|X|^{1/2}(C|X|^{1/2}-1)\geq \frac{1}{2}C^2|Y||X|$. This means that on average a pair of distinct vertices in $X$ have at least $\frac{1}{2}C^2\frac{|Y|}{|X|}$ common neighbours. Hence, the number of $K_{2,2}$'s in $H$ is at least $\binom{|X|}{2}\binom{\frac{1}{2}C^2\frac{|Y|}{|X|}}{2}\geq \frac{1}{1000}C^4|Y|^2$. Since $|X|<|Y|$ and $C$ is sufficiently large, Lemma \ref{lemma:blownuptree} shows that $H$ contains a copy of $T[2]$ with the two images of $v$ embedded in $X$.
\end{proof}

We can now prove Lemma \ref{lemma:3heavy}.

\begin{proof}[Proof of Lemma \ref{lemma:3heavy}]
	Let $C$ be sufficiently large and suppose that the number of triples $(z',w_1,w_2)$ of distinct vertices with $z'\in R$, $w_1,w_2\in N(z',z_1,z_2)$ is more than $C\ell^2$. By Lemma \ref{lemma:count cherries} with $X=N(z_1,z_2)$, $Y=R$ and $H=G[X\cup Y]$, there exist distinct vertices $u,u',w_1,w_1',w_2,w_2',\dots,w_t,w_t'\in N(z_1,z_2)$ and $v_1,v_1',v_2,v_2',\dots,v_t,v_t'\in R$ such that $v_i,v_i'\in N(u,u',w_i,w_i')$ for every $1\leq i\leq t$. Then the vertices $z_1,z_2,u,u',w_1,w_1',w_2,w_2',\dots,w_t,w_t',v_1,v_1',v_2,v_2',\dots,v_t,v_t'$ together form a copy of $\blownuptheta$ (see Figure \ref{fig:3heavy}), which is a contradiction.
	\begin{figure}
		\centering
		\begin{tikzpicture}[scale=0.7]
		\draw[fill=black](6,6)circle(5pt);
		\draw[fill=black](7,6)circle(5pt);
		
		\draw[fill=black](3,6)circle(5pt);
		\draw[fill=black](2,6)circle(5pt);
		\draw[fill=black](-2,6)circle(5pt);
		\draw[fill=black](-5,6)circle(5pt);
		\draw[fill=black](-1,6)circle(5pt);
		\draw[fill=black](-4,6)circle(5pt);
		
		\draw[fill=black](-6,12)circle(5pt);
		\draw[fill=black](-5,12)circle(5pt);
		\draw[fill=black](-2,12)circle(5pt);
		\draw[fill=black](-1,12)circle(5pt);
		\draw[fill=black](1,12)circle(5pt);
		\draw[fill=black](2,12)circle(5pt);
		\draw[fill=black](5,12)circle(5pt);
		\draw[fill=black](6,12)circle(5pt);
		
		\draw[fill=black](3.1,12)circle(1pt);
		\draw[fill=black](3.4,12)circle(1pt);
		\draw[fill=black](3.7,12)circle(1pt);		
		\draw[fill=black](0.1,6)circle(1pt);
		\draw[fill=black](0.4,6)circle(1pt);
		\draw[fill=black](0.7,6)circle(1pt);
					
		\draw[thick](-6,12)--(-5,6)--(-2,12)--(6,6);
		\draw[thick](-6,12)--(-4,6)--(-2,12)--(7,6);
		\draw[thick](-5,12)--(-5,6)--(-1,12)--(6,6);
		\draw[thick](-5,12)--(-4,6)--(-1,12)--(7,6);
		
		\draw[thick](-6,12)--(-2,6)--(1,12)--(6,6);
		\draw[thick](-6,12)--(-1,6)--(1,12)--(7,6);
		\draw[thick](-5,12)--(-2,6)--(2,12)--(6,6);
		\draw[thick](-5,12)--(-1,6)--(2,12)--(7,6);
		
		\draw[thick](-6,12)--(2,6)--(5,12)--(6,6);
		\draw[thick](-6,12)--(3,6)--(5,12)--(7,6);
		\draw[thick](-5,12)--(2,6)--(6,12)--(6,6);
		\draw[thick](-5,12)--(3,6)--(6,12)--(7,6);
		
		\node at (6,5.5) {$z_1$};
		\node at (7,5.5) {$z_2$};
		
		\node at (2,5.5) {$v_t$};
		\node at (3,5.5) {$v_t'$};
		\node at (-2,5.5) {$v_2$};
		\node at (-1,5.5) {$v_2'$};
		\node at (-5,5.5) {$v_1$};
		\node at (-4,5.5) {$v_1'$};

		\node at (-2,12.5) {$w_1$};
		\node at (-1,12.5) {$w_1'$};
		\node at (1,12.5) {$w_2$};
		\node at (2,12.5) {$w_2'$};
		\node at (5,12.5) {$w_t$};
		\node at (6,12.5) {$w_t'$};
		
		\node at (-6,12.4) {$u$};
		\node at (-5,12.5) {$u'$};
		
		\draw[rotate around={90:(-1,6)}] (-1,6) ellipse (1 and 6);
		\draw[rotate around={90:(0,12)}] (0,12) ellipse (1 and 8);
		\node at (-10,12)  {$N(z_1,z_2)$};
		\node at (-8,6)  {$R$};
		\end{tikzpicture}		
		\caption{$\blownuptheta$ in the proof of Lemma \ref{lemma:3heavy}}
		\label{fig:3heavy}
	\end{figure}
\end{proof}

\subsection{Finding $P_3[2]$'s with prescribed extra properties} \label{subsec:prescribeproperties}

This is the subsection where most of the properties prescribed for our $P_3[2]$'s will be established. Roughly speaking, the next lemma will be used to find $P_3[2]$'s with the property that $d(x_1,x_2,z_1)$ and $d(x_1,x_2,z_2)$ are not too large (as before, we are referring to Figure \ref{fig:blownuppath} here).

\begin{lemma} \label{lemma:3light}
	Let $G$ be a graph on $n$ vertices with minimum degree $\d=\omega(n^{2/3})$. Let $S\subset V(G)$ have size $s\geq n^{1/3}$. Then there exists some $\lambda=\omega(1)$ such that the number of vertices $v\in V(G)$ with $\frac{\lambda}{2} \frac{s}{n^{1/3}}< d_S(v)\leq \lambda \frac{s}{n^{1/3}}$ is at least $c\delta n^{1/3}\lambda^{-11/10}$, where $c=(\sum_{i\geq 0}2^{-i/10})^{-1}$.
\end{lemma}

\begin{proof}
	Define $U_0=\{v\in V(G): d_S(v)\leq  \frac{s}{n^{1/3}}\}$, and for every positive integer $i$, let $U_i=\{v\in V(G): \frac{s}{n^{1/3}}2^{i-1}< d_S(v)\leq  \frac{s}{n^{1/3}}2^{i}\}$.
	
	Now we double count the number of edges between $S$ and $V(G)$ (viewed as a bipartite graph). On the one hand, every $y\in S$ has at least $\d$ neighbours in $V(G)$. On the other hand, any $v\in U_i$ has at most $\frac{s}{n^{1/3}}2^{i}$ neighbours in $S$. Thus,
	$$\sum_{i\geq 0} |U_i|\frac{s}{n^{1/3}}2^{i}\geq s \delta,$$
	so
	$$\sum_{i\geq 0} |U_i|2^{i}\geq \delta n^{1/3}.$$
	It is easy to see that then there exists some $i$ such that $|U_i|\geq 2^{-\frac{11i}{10}}c\delta n^{1/3}$. Since $|U_i|\leq n$, we have $i=\omega(1)$. So we may take $\lambda=2^{i}$.
\end{proof}

The next lemma lists almost all properties that we require about the vertices $x_1,x_2,$ $y_1,y_2,$ $z_1,z_2,$ $w_1, w_2$. The one additional property that we will need is that $d(y_1,y_2,w_1,w_2)<6t$.

\begin{lemma} \label{lemma:goodblownuppaths}
	Let $K$ be a constant and let $G$ be a $K$-almost-regular, $\blownuptheta$-free graph on $n$ vertices with minimum degree $\d=\omega(n^{2/3})$. Then there exist distinct vertices $x_1,x_2$ in $G$ and a set $S\subset N(x_1,x_2)$ of size at least $n^{1/3}$ as follows. Writing $s=|S|$, there exist $\lambda=\omega(1)$, $\mu=\omega(1)$, and a collection $\mathcal{Q}$ of $\Omega(\frac{s^2n^{2/3}\lambda^{27/10}}{\mu^{11/10}})$ tuples $(y_1,y_2,z_1,z_2)\in V(G)^4$ satisfying the following properties.
	
	\begin{enumerate}
		\item $y_1,y_2\in S$ and $y_iz_j$ are edges for every $i,j$. \label{prop:edges}
		\item $x_1,x_2,y_1,y_2,z_1,z_2$ are distinct. \label{prop:distinct}
		\item $d_S(z_1),d_S(z_2)\leq \lambda\frac{s}{n^{1/3}}$. \label{prop:triple codegree}
		\item $d(x_1,x_2,z_1,z_2)<6t$. \label{prop:quadruple codegree}
		\item $\mu n^{1/3}\leq d(z_1,z_2)\leq 2\mu n^{1/3}$. \label{prop:pair codegree}
	\end{enumerate}
\end{lemma}

\begin{proof}
	Let $c=(\sum_{i\geq 0}2^{-i/10})^{-1}$ as in Lemma \ref{lemma:3light}. For every $R\subset V(G)$ of size at least $n^{1/3}$, choose $\lambda(R)=\omega(1)$ such that the number of vertices $v$ with $\frac{\lambda(R)}{2} \frac{|R|}{n^{1/3}}< d_R(v)\leq \lambda(R) \frac{|R|}{n^{1/3}}$ is at least $c\delta n^{1/3}\lambda(R)^{-11/10}$. By Lemma \ref{lemma:3light}, such a choice exists. Since $G$ has minimum degree $\omega(n^{2/3})$, it is easy to see that if $n$ is sufficiently large, then there exist distinct $u,v\in V(G)$ with $d(u,v)\geq n^{1/3}$.
	Choose distinct $x_1,x_2\in V(G)$ and $S\subset N(x_1,x_2)$ such that $|S|\geq n^{1/3}$ and $\lambda(S)$ is minimal among these choices. Let $\lambda=\lambda(S)$. It remains to find $\mu$ and enough tuples $(y_1,y_2,z_1,z_2)$ with properties \ref{prop:edges}-\ref{prop:pair codegree}. This is done in two main steps.
	
	\medskip
	
	\textbf{Step 1.} We find $\Omega(s^2 n^{2/3} \lambda^{9/5})$ tuples $(y_1,y_2,z_1,z_2)$ satisfying properties~\ref{prop:edges},~\ref{prop:distinct} and~\ref{prop:triple codegree}.
	
	\medskip
	
	Let $U=\{v\in V(G)\setminus \{x_1,x_2\}: \frac{\lambda}{2} \frac{s}{n^{1/3}}< d_S(v)\leq  \lambda \frac{s}{n^{1/3}}\}$. Then $|U|\geq c\delta n^{1/3}\lambda^{-11/10}-2\geq n \lambda^{-11/10}$ for $n$ sufficiently large.
	
	Clearly, the number of triples $(y_1,y_2,z)$ with $y_1,y_2\in S$ distinct, $z\in U$ and $y_1z,y_2z\in E(G)$ is at least $|U|(\frac{\lambda}{2}\frac{s}{n^{1/3}})(\frac{\lambda}{2}\frac{s}{n^{1/3}}-1)=\Omega(s^2 n^{1/3} \lambda^{9/10})$. Hence, on average, for a pair $y_1,y_2\in S$ there are $\Omega(n^{1/3}\lambda^{9/10})$ vertices $z\in N(y_1,y_2)\cap U$. By convexity, on average, for a pair $y_1,y_2\in S$ there are $\Omega(n^{2/3}\lambda^{9/5})$ pairs of distinct vertices $z_1,z_2\in N(y_1,y_2)\cap U$. Since any $z\in U$ has $d_S(z)\leq \lambda\frac{s}{n^{1/3}}$, this completes Step~1.
	
	\medskip
	
	\textbf{Step 2.} We find $\Omega(s^2 n \lambda^{27/10})$ tuples $(y_1,y_2,z_1,z_2,w)$ satisfying properties \ref{prop:edges}, \ref{prop:distinct}, \ref{prop:triple codegree} and \ref{prop:quadruple codegree} with the additional properties that $d(z_1,z_2)\geq n^{1/3}\lambda^{4/5}$ and $z_1w,z_2w\in E(G)$.
	
	\medskip
	
	For $y_1,y_2\in S$, let $N(y_1,y_2)^*=\{v\in N(y_1,y_2)\setminus \{x_1,x_2\}: d_S(v)\leq \lambda\frac{s}{n^{1/3}}\}$. The conclusion of Step 1 implies that
	$$\sum_{y_1,y_2\in S \text{ distinct }} |N(y_1,y_2)^*|^2=\Omega(s^2 n^{2/3} \lambda^{9/5}).$$
	Hence,
	\begin{equation}
		\sum_{\substack{y_1,y_2\in S \text{ distinct } \\ |N(y_1,y_2)^*|\geq n^{1/3}}} |N(y_1,y_2)^*|^2=\Omega(s^2 n^{2/3} \lambda^{9/5}). \label{eqn:step1conclusion}
	\end{equation}
	We now prove that for any distinct $y_1,y_2\in S$ with $|N(y_1,y_2)^*|\geq n^{1/3}$, the number of triples $(z_1,z_2,w)$ of distinct vertices with $z_1,z_2\in N(y_1,y_2)^*$, $d(x_1,x_2,z_1,z_2)<6t$, $d(z_1,z_2)\geq n^{1/3}\lambda^{4/5}$ and $w\in N(z_1,z_2)$ is $\Omega(|N(y_1,y_2)^*|^2 n^{1/3}\lambda^{9/10})$. Using equation (\ref{eqn:step1conclusion}), this would complete Step 2.
	
	Let some distinct $y_1,y_2\in S$ have $|N(y_1,y_2)^*|\geq n^{1/3}$. Let $R=N(y_1,y_2)^*$. By definition, the number of vertices $v$ with $d_R(v)>\frac{\lambda(R)}{2}\frac{|R|}{n^{1/3}}$ is at least $c \delta n^{1/3}\lambda(R)^{-11/10}$. Thus, the number of triples of distinct vertices $(z_1,z_2,w)$ with $z_1,z_2\in R$ and $w\in N(z_1,z_2)$ is $\Omega(|R|^2 \delta n^{-1/3} \lambda(R)^{9/10})\geq \Omega(|R|^2 \delta n^{-1/3} \lambda^{9/10})$.
	By Lemma \ref{lemma:4heavy}, the number of pairs of distinct vertices $(z_1,z_2)\in R^2$ with $d(x_1,x_2,z_1,z_2)\geq 6t$ is at most $4t|R|$. Since $G$ has maximum degree at most $K\delta$, the number of triples $(z_1,z_2,w)$ involving such pairs $(z_1,z_2)$ is at most $4t|R|K\d$. Note that $|R|\geq n^{1/3}$ and $\lambda =\omega(1)$, so $4t|R|K\d=o(|R|^2 \delta n^{-1/3} \lambda^{9/10})$. Moreover, the number of triples $(z_1,z_2,w)$ with $z_1,z_2\in R$, $w\in N(z_1,z_2)$ and $d(z_1,z_2)\leq n^{1/3}\lambda^{4/5}$ is clearly at most $|R|^2 n^{1/3} \lambda^{4/5}$, which is again $o(|R|^2 \delta n^{-1/3} \lambda^{9/10})$. Thus, the number of triples $(z_1,z_2,w)$ of distinct vertices with $z_1,z_2\in R$, $d(x_1,x_2,z_1,z_2)<6t$, $d(z_1,z_2)\geq n^{1/3} \lambda^{4/5}$ and $w\in N(z_1,z_2)$ is $\Omega(|R|^2 \delta n^{-1/3} \lambda^{9/10})$. This is $\Omega(|N(y_1,y_2)^*|^2 n^{1/3} \lambda^{9/10})$, as claimed.
	
	\medskip
	
	Since $\sum_{j=1}^{\infty} 2^{-j/10}$ is bounded, the conclusion of Step 2 implies that there exists some positive integer $j$ such that there exist $\Omega(\frac{s^2 n \lambda^{27/10}}{2^{j/10}})$ tuples $(y_1,y_2,z_1,z_2,w)$ satisfying properties \ref{prop:edges}, \ref{prop:distinct}, \ref{prop:triple codegree} and \ref{prop:quadruple codegree} with the additional properties $n^{1/3} \lambda^{4/5} 2^{j-1}\leq d(z_1,z_2)<n^{1/3} \lambda^{4/5} 2^j$ and $z_1w,z_2w\in E(G)$. Take $\mu=\lambda^{4/5} 2^{j-1}=\omega(1)$. Then the number of tuples $(y_1,y_2,z_1,z_2)$ satisfying properties \ref{prop:edges}-\ref{prop:pair codegree} is $\Omega(\frac{s^2 n \lambda^{27/10}}{2^{j/10}n^{1/3}\lambda^{4/5}2^j})\geq \Omega(\frac{s^2 n^{2/3} \lambda^{27/10}}{\mu^{11/10}})$.
\end{proof}

\subsection{Counting the number of pairs of $P_3[2]$'s which share an internal vertex} \label{subsec:countbadpairsofP3}

As explained in the outline of the proof, we now argue that if $G$ does not contain $\blownuptheta$ as a subgraph, then there are plenty of copies of the configuration depicted in Figure~\ref{fig:degenerate blownupc6}. The next lemma makes this precise.

\begin{lemma} \label{lemma:number of overlapping pairs}
    Let $G$ be a $\blownuptheta$-free bipartite graph on $n$ vertices. Let $x_1,x_2$ be distinct vertices, let $S\subset N(x_1,x_2)$ be a set of size $s\geq n^{1/3}$ and let $\lambda,\mu=\omega(1)$. Assume that there is a set $\mathcal{Q}$ of $q=\Omega(\frac{s^2n^{2/3}\lambda^{27/10}}{\mu^{11/10}})$ tuples $(y_1,y_2,z_1,z_2)\in V(G)^4$ satisfying the five properties in Lemma \ref{lemma:goodblownuppaths}. Then there exists a set $\mathcal{A}$ of $\Omega(\mu^4 n^{-2/3}q^2)$ tuples $(y_1,y_2,z_1,z_2,w_1,w_2,y_2',z_1',z_2')\in V(G)^9$ for which the following properties hold. 
    
    \begin{enumerate}[label=(\alph*)]
        \item $(y_1,y_2,z_1,z_2),(y_1,y_2',z_1',z_2')\in \mathcal{Q}$.
        \item $w_1,w_2\in N(z_1,z_2,z_1',z_2')$.
        \item $d(y_1,y_2,w_1,w_2),d(y_1,y_2',w_1,w_2)<6t$.
        \item $x_1,x_2,y_1,y_2,z_1,z_2,w_1,w_2,y_2',z_1',z_2'$ are all distinct.
    \end{enumerate}
    \begin{comment}
    \begin{figure}
		\centering
		\begin{tikzpicture}[scale=0.6]
		\draw[fill=black](0,0)circle(5pt);
		\draw[fill=black](0,2)circle(5pt);		
		\draw[fill=black](6,-1)circle(5pt);
		\draw[fill=black](6,1)circle(5pt);
		\draw[fill=black](6,3)circle(5pt);
		\draw[fill=black](18,0)circle(5pt);
		\draw[fill=black](18,2)circle(5pt);
		\draw[fill=black](12,4)circle(5pt);
		\draw[fill=black](12,2)circle(5pt);
		\draw[fill=black](12,0)circle(5pt);
		\draw[fill=black](12,-2)circle(5pt);
		
		\draw[thick](0,2)--(6,3)--(12,4)--(18,2);
		\draw[thick](0,2)--(6,1)--(12,4)--(18,0);
		\draw[thick](0,0)--(6,3)--(12,2)--(18,2);
		\draw[thick](0,0)--(6,1)--(12,2)--(18,0);
		
		\draw[thick](0,0)--(6,-1)--(12,0)--(18,0);
		\draw[thick](0,2)--(6,-1)--(12,-2)--(18,0);
		\draw[thick](6,1)--(12,0)--(18,2);
		\draw[thick](6,1)--(12,-2)--(18,2);
		
		\node at (0,-0.7) {$x_1$};
		\node at (0,1.3) {$x_2$};
		\node at (6,-1.7) {$y_2'$};
		\node at (6,2.3) {$y_2$};
		\node at (6,0.3)  {$y_1=y_1'$};
		\node at (12,3.3)  {$z_2$};
		\node at (12,-2.7)  {$z_2'$};
		\node at (12,-0.7)  {$z_1'$};
		\node at (12,1.3) {$z_1$};
		\node at (18,-0.7) {$w_1$};
		\node at (18,1.3) {$w_2$};
		\end{tikzpicture}
		\caption{An element of $\mathcal{A}$}
		\label{fig:gluedblownuppaths}
	\end{figure}
	\end{comment}
\end{lemma}

\begin{proof}
    By property \ref{prop:pair codegree} from Lemma \ref{lemma:goodblownuppaths} and Lemma \ref{lemma:4heavy}, any $(y_1,y_2,z_1,z_2)\in \mathcal{Q}$ can be extended in $\Theta(\mu^2 n^{2/3})$ ways to a tuple $(y_1,y_2,z_1,z_2,w_1,w_2)$ of vertices with the additional properties that $w_1$ and $w_2$ are distinct elements of $N(z_1,z_2)\setminus \{x_1,x_2,y_1,y_2\}$ and $d(y_1,y_2,w_1,w_2)<6t$. Let $\mathcal{R}$ be the set of all tuples obtained this way and let $r=|\mathcal{R}|$. Note that $r=\Theta(q\mu^2 n^{2/3})$, so $r=\omega((\frac{\lambda s}{n^{1/3}})^2 n^2)$. Thus, on average a pair $(w_1,w_2)$ of distinct vertices can be extended in $\omega((\frac{\lambda s}{n^{1/3}})^2)$ ways to a tuple $(y_1,y_2,z_1,z_2,w_1,w_2)\in \mathcal{R}$.
	
	Fix a pair $(w_1,w_2)\in V(G)^2$ and assume it can be extended to $h=\omega((\frac{\lambda s}{n^{1/3}})^2)$ such tuples. Find a maximal set of disjoint tuples $(y_1^1,y_2^1,z_1^1,z_2^1),(y_1^2,y_2^2,z_1^2,z_2^2),\dots,(y_1^k,y_2^k,z_1^k,z_2^k)$ such that $(y_1^{i},y_2^{i},z_1^{i},z_2^{i},w_1,w_2)\in \mathcal{R}$ for every $1\leq i\leq k$. Since $G$ is $\blownuptheta$-free, we have $k<t$. Now for any $y_1,y_2,z_1,z_2$ with $(y_1,y_2,z_1,z_2,w_1,w_2)\in \mathcal{R}$, we have $\{y_1,y_2,z_1,z_2\}\cap \{y_1^1,y_2^1,z_1^1,z_2^1,y_1^2,y_2^2,z_1^2,z_2^2,\dots,y_1^k,y_2^k,z_1^k,z_2^k\}\neq \emptyset$. By the pigeon hole principle, there exists some $v\in \{y_1^1,y_2^1,z_1^1,z_2^1,y_1^2,y_2^2,z_1^2,z_2^2,\dots,y_1^k,y_2^k,z_1^k,z_2^k\}$ such that at least one of the following holds.
	\begin{enumerate}[label=(\roman*)]
		\item There are at least $\frac{h}{16k}$ tuples $(y_1,y_2,z_1,z_2,w_1,w_2)\in \mathcal{R}$ with $y_1=v$.
		\item There are at least $\frac{h}{16k}$ tuples $(y_1,y_2,z_1,z_2,w_1,w_2)\in \mathcal{R}$ with $y_2=v$.
		\item There are at least $\frac{h}{16k}$ tuples $(y_1,y_2,z_1,z_2,w_1,w_2)\in \mathcal{R}$ with $z_1=v$.
		\item There are at least $\frac{h}{16k}$ tuples $(y_1,y_2,z_1,z_2,w_1,w_2)\in \mathcal{R}$ with $z_2=v$.
	\end{enumerate}
	If $(y_1,y_2,z_1,z_2,w_1,w_2)\in \mathcal{R}$, then $y_1,y_2\in S$ (by the definition of $\mathcal{Q}$ and property \ref{prop:edges} in Lemma \ref{lemma:goodblownuppaths}), $d_S(z_1)\leq \frac{\lambda s}{n^{1/3}}$ (by property \ref{prop:triple codegree}), and $d(y_1,y_2,w_1,w_2)<6t$ (by the definition of $\mathcal{R}$). Thus, there are at most $(\frac{\lambda s}{n^{1/3}})^2\cdot 6t$ ways to extend a fixed choice of $z_1,w_1,w_2$ to get $(y_1,y_2,z_1,z_2,w_1,w_2)\in \mathcal{R}$. In particular, (since $h=\omega((\frac{\lambda s}{n^{1/3}})^2)$), case (iii) is impossible for $n$ sufficiently large. Similarly, case (iv) is impossible. Thus, either case (i) or case (ii) holds.
	
	Assume, without loss of generality, that case (i) holds. Since $d(y_1,y_2,w_1,w_2)<6t$, for any $u\in V(G)$ there are at most $(6t)^2$ tuples $(y_1,y_2,z_1,z_2,w_1,w_2)\in \mathcal{R}$ with $y_1=v,y_2=u$. Moreover, for any $u\in V(G)$ there are at most $\frac{\lambda s}{n^{1/3}}\cdot 6t$ tuples $(y_1,y_2,z_1,z_2,w_1,w_2)\in \mathcal{R}$ with $y_1=v,z_1=u$, and there are at most $\frac{\lambda s}{n^{1/3}}\cdot 6t$ tuples $(y_1,y_2,z_1,z_2,w_1,w_2)\in \mathcal{R}$ with $y_1=v,z_2=u$.  Hence, almost all pairs from our at least $\frac{h}{16k}$ tuples $(y_1,y_2,z_1,z_2,w_1,w_2)\in \mathcal{R}$ with $y_1=v$ are disjoint apart from $y_1$, $w_1$ and $w_2$. Thus, for our fixed $w_1,w_2$, there are $\Omega(h^2)$ pairs $(y_1,y_2,z_1,z_2,w_1,w_2),(y_1',y_2',z_1',z_2',w_1,w_2)\in \mathcal{R}$ with $y_1=y_1'$ but $\{y_2,z_1,z_2\}\cap \{y_2',z_1',z_2'\}=\emptyset$.

	Summing over all pairs $(w_1,w_2)$ and noting the symmetry of cases (i) and (ii) above, we get $\Omega(n^2\cdot (\frac{r}{n^2})^2)=\Omega(\frac{r^2}{n^2})$ pairs $(y_1,y_2,z_1,z_2,w_1,w_2),(y_1',y_2',z_1',z_2',w_1,w_2)\in \mathcal{R}$ with $y_1=y_1'$ but $\{y_2,z_1,z_2\}\cap \{y_2',z_1',z_2'\}=\emptyset$.
	Let $\mathcal{A}$ be the set of all tuples $(y_1,y_2,z_1,z_2,w_1,w_2,y_2',z_1',z_2')$ for which $(y_1,y_2,z_1,z_2,w_1,w_2),(y_1,y_2',z_1',z_2',w_1,w_2)\in \mathcal{R}$ and $\{y_2,z_1,z_2\}\cap \{y_2',z_1',z_2'\}=\emptyset$. Then $|\mathcal{A}|=\Omega(\frac{r^2}{n^2})=\Omega(\mu^4 n^{-2/3}q^2)$.
\end{proof}

We will now force a contradiction to the $\blownuptheta$-freeness assumption by showing that
\begin{itemize}
    \item there are ``few" tuples in $\mathcal{A}$ with ``large" $d(z_1,z_1',z_2')$ or ``large" $d(z_2,z_1',z_2')$, and
    \item there are ``few" tuples in $\mathcal{A}$ with ``small" $d(z_1,z_1',z_2')$ and ``small" $d(z_2,z_1',z_2')$.
\end{itemize}

The next lemma establishes the former statement.

\begin{lemma} \label{lemma:number with large 3-degree}
    In the setting of Lemma \ref{lemma:number of overlapping pairs}, the number of tuples $(y_1,y_2,z_1,z_2,w_1,w_2,y_2',z_1',z_2')\in \mathcal{A}$ with $d(z_1,z_1',z_2')\geq \mu^{7/5}n^{1/6}$ is $o(|\mathcal{A}|)$.\\
    Moreover, the number of tuples $(y_1,y_2,z_1,z_2,w_1,w_2,y_2',z_1',z_2')\in \mathcal{A}$ with $d(z_2,z_1',z_2')\geq \mu^{7/5}n^{1/6}$ is also $o(|\mathcal{A}|)$.
\end{lemma}

\begin{proof}
    By symmetry, it suffices to prove the first statement.

    Let $(y_1,y_2,z_1,z_2,w_1,w_2,y_2',z_1',z_2')\in \mathcal{A}$ with $d(z_1,z_1',z_2')\geq \mu^{7/5}n^{1/6}$. We bound the number of possibilities as follows. First note that $(y_1,y_2',z_1',z_2')\in \mathcal{Q}$, so there are at most $q$ choices for these vertices. For any such choice, $\mu n^{1/3}\leq d(z_1',z_2')\leq 2\mu n^{1/3}$ holds by property \ref{prop:pair codegree} from Lemma \ref{lemma:goodblownuppaths}. Since $\mu^{7/5}n^{1/6}=\omega((2\mu n^{1/3})^{1/2})$, Lemma~\ref{lemma:3heavy} (applied with $z_1'$ and $z_2'$ in place of $z_1$ and $z_2$) implies that there are $O((\mu n^{1/3})^2)$ choices for $(z_1,w_1,w_2)$. Moreover, there are at most $d_S(z_1)\leq \lambda\frac{s}{n^{1/3}}$ choices for $y_2$. Finally, there are at most $d(y_1,y_2,w_1,w_2)<6t$ choices for $z_2$. Altogether, we find that there are $O(q\cdot (\mu n^{1/3})^2\cdot \lambda\frac{s}{n^{1/3}}\cdot 6t)$ elements of $\mathcal{A}$ with $d(z_1,z_1',z_2')\geq \mu^{7/5}n^{1/6}$. Since $|\mathcal{A}|=\Omega(\mu^4n^{-2/3}q^2)$, $q=\Omega(\frac{s^2 n^{2/3} \lambda^{27/10}}{\mu^{11/10}})$, $s\geq n^{1/3}$, $\lambda=\omega(1)$ and $\mu=\omega(1)$, we have
    $$q\cdot (\mu n^{1/3})^2\cdot \lambda\frac{s}{n^{1/3}}\cdot 6t=o(|\mathcal{A}|),$$
    completing the proof.
\end{proof}

We now turn to counting the number of tuples in $\mathcal{A}$ with small $d(z_1,z_1',z_2')$ and small $d(z_2,z_1',z_2')$. In the setting of Lemma \ref{lemma:number of overlapping pairs}, let $\mathcal{A}'$ be the set of $(y_1,y_2,z_1,z_2,w_1,w_2,y_2',z_1',z_2')\in \mathcal{A}$ with $d(z_1,z_1',z_2')< \mu^{7/5}n^{1/6}$ and $d(z_2,z_1',z_2')< \mu^{7/5}n^{1/6}$.

\begin{lemma}\label{lemma:final}
    If $n$ is sufficiently large, then any $y_1,y_2,y_2',z_1',z_2'$ can be extended to at most $g:=10^8t^3\mu^{14/5}n^{1/3}$ elements of $\mathcal{A}'$.
\end{lemma}
    
\begin{proof}

    Suppose otherwise. Take a maximal set of disjoint tuples $(z_1^1,z_2^1,w_1^1,w_2^1),(z_1^2,z_2^2,w_1^2,w_2^2),\dots,(z_1^k,z_2^k,w_1^k,w_2^k)$ such that $(y_1,y_2,z_1^{i},z_2^{i},w_1^{i},w_2^{i},y_2',z_1',z_2')\in \mathcal{A}'$ for every $1\leq i\leq k$. Since $G$ is $\blownuptheta$-free, we have $k<t$. Now for any $z_1,z_2,w_1,w_2$ with $(y_1,y_2,z_1,z_2,w_1,w_2,y_2',z_1',z_2')\in \mathcal{A}'$, we have $\{z_1,z_2,w_1,w_2\}\cap \{z_1^1,z_2^1,w_1^1,w_2^1,z_1^2,z_2^2,w_1^2,w_2^2,\dots,z_1^k,z_2^k,w_1^k,w_2^k\}\neq \emptyset$. By the pigeonhole principle, there exists some $v\in \{z_1^1,z_2^1,w_1^1,w_2^1,z_1^2,z_2^2,w_1^2,w_2^2,\dots,z_1^k,z_2^k,w_1^k,w_2^k\}$ such that at least one of the following holds.
    \begin{enumerate}[label=(\roman*)]
		\item There are at least $\frac{g}{16k}$ ways to extend $y_1,y_2,y_2',z_1',z_2'$ to an element of $\mathcal{A}'$ with $z_1=v$.
		\item There are at least $\frac{g}{16k}$ ways to extend $y_1,y_2,y_2',z_1',z_2'$ to an element of $\mathcal{A}'$ with $z_2=v$.
		\item There are at least $\frac{g}{16k}$ ways to extend $y_1,y_2,y_2',z_1',z_2'$ to an element of $\mathcal{A}'$ with $w_1=v$.
		\item There are at least $\frac{g}{16k}$ ways to extend $y_1,y_2,y_2',z_1',z_2'$ to an element of $\mathcal{A}'$ with $w_2=v$.
	\end{enumerate}
	Assume first that case (i) holds. Recall that the definition of $\mathcal{A}'$ requires that $d(z_1,z_1',z_2')<\mu^{7/5}n^{1/6}$. This means that there are at most $(\mu^{7/5}n^{1/6})^2$ ways to choose $w_1$ and $w_2$, given $y_1,y_2,y_2',z_1,z_1',z_2'$. Furthermore, we must have $d(y_1,y_2,w_1,w_2)<6t$, so, given fixed choices for $w_1$ and $w_2$ there are at most $6t$ ways to choose $z_2$. Altogether this allows at most $(\mu^{7/5}n^{1/6})^2\cdot 6t$ ways to extend $y_1,y_2,y_2',z_1,z_1',z_2'$ to an element of $\mathcal{A'}$. This is less than $\frac{g}{16k}$, a contradiction.
	Similarly case (ii) leads to a contradiction.
	
	Assume now that case (iii) holds. In any element of $\mathcal{A}$, we have (by property \ref{prop:pair codegree} from Lemma \ref{lemma:goodblownuppaths}) $d(z_1',z_2')\leq 2\mu n^{1/3}$ and $w_2\in N(z_1',z_2')$. Moreover, $z_1,z_2\in N(y_1,y_2,w_1,w_2)$ and $d(y_1,y_2,w_1,w_2)<6t$, so there are at most $2\mu n^{1/3}\cdot (6t)^2$ ways to extend $y_1,y_2,w_1,y_2',z_1',z_2'$ to an element of $\mathcal{A}$. Since $\mu=\omega(1)$, this is again less than $\frac{g}{16k}$ for sufficiently large $n$, a contradiction. Similarly case (iv) leads to a contradiction.
	\end{proof}

	\begin{lemma} \label{lemma:number with small 3-degree}
    In the setting of Lemma \ref{lemma:number of overlapping pairs}, there are $o(|\mathcal{A}|)$  tuples $(y_1,y_2,z_1,z_2,w_1,w_2,y_2',z_1',z_2')\in \mathcal{A}$ for which $d(z_1,z_1',z_2')< \mu^{7/5}n^{1/6}$ and $d(z_2,z_1',z_2')< \mu^{7/5}n^{1/6}$.
\end{lemma}

\begin{proof}	Note first that for any $(y_1,y_2,z_1,z_2,w_1,w_2,y_2',z_1',z_2')\in \mathcal{A}$, we have $(y_1,y_2',z_1',z_2')\in \mathcal{Q}$ and $y_2\in S$, so there are at most $qs$ choices for $y_1,y_2,y_2',z_1',z_2'$. Combining this with Lemma \ref{lemma:final}, we find that $|\mathcal{A}'|\leq 10^8t^3\mu^{14/5}n^{1/3}qs$. Since $|\mathcal{A}|=\Omega(\mu^4n^{-2/3}q^2)$, $q=\Omega(\frac{s^2 n^{2/3} \lambda^{27/10}}{\mu^{11/10}})$, $s\geq n^{1/3}$, $\lambda=\omega(1)$ and $\mu=\omega(1)$, we obtain $|\mathcal{A}'|=o(|\mathcal{A}|)$.
\end{proof}

Observe that Theorem \ref{thm:turanthetareduced} follows by combining Lemma~\ref{lemma:goodblownuppaths}, Lemma~\ref{lemma:number of overlapping pairs}, Lemma~\ref{lemma:number with large 3-degree} and Lemma~\ref{lemma:number with small 3-degree}.

\section{Concluding remarks}

Conjecture \ref{con:blownupcycles} is still open is most cases and appears to be quite hard. Even proving the conjectured upper bound for $\ex(n,C_6[r])$ (with general $r$) or $\ex(n,C_{2k}[2])$ (with general $k$) seems to require some new ideas.

Another interesting problem is to give good lower bounds for $\ex(n,C_6[2])$ and to find a small value of $t$ for which $\ex(n,\theta_{3,t}[2])=\Theta(n^{5/3})$ already holds. 

\vspace{4mm}

\noindent \textbf{Acknowledgement.} We are grateful to Craig Timmons for discussions on constructions using finite fields and for his generous help with computations.

\end{document}